\documentclass[11pt]{amsart}
\usepackage{amsmath,amsthm, amsfonts,amssymb}
\date{}
\newtheorem{theorem}{Theorem}%[section]
\newtheorem{lemma}{Lemma}%[section]

\newtheorem{corollary}{Corollary}%[section]
\newtheorem{proposition}{Proposition}%[section]

\newtheorem{problem}{Problem}

\theoremstyle{definition}
\newcommand{\sep}{, }
\newcommand{\Bd}{\mathrm{Bd}}

\newcommand{\Dis}{\mathcal{D}}
\newcommand{\uDis}{u\mathcal{D}}

\newcommand{\diam}{\mathrm{diam}}
\newcommand{\IN}{\mathbb{N}}

\newcommand{\IR}{\mathbb{R}}
\newcommand{\IQ}{\mathbb{Q}}
\newcommand{\w}{\omega}
\newcommand{\U}{\mathcal{U}}

\newcommand{\Hip}{\mathcal{H}}
\newcommand{\e}{\varepsilon}
\newcommand{\cl}{\mathrm{cl}}
\newcommand{\dist}{\mathrm{dist}}

\begin{document}
%\begin{frontmatter}
\title[Characterizing ANR-property in hyperspaces]
{Characterizing metric spaces whose hyperspaces are absolute
neighborhood retracts}
\author{T.Banakh and R.Voytsitskyy}
\address{Instytut Matematyki, Akademia \'Swi\c etokrzyska, Kielce,
Poland and\newline Department of Mathematics, Ivan Franko Lviv
National University, Lviv, Ukraine} \email{tbanakh@franko.lviv.ua}
\begin{abstract}
We characterize metric spaces $X$ whose
hyperspaces $2^X$ or $\Bd(X)$ of non-empty closed (bounded)
subsets, endowed with the Hausdorff metric, are absolute
[neighborhood] retracts.
\end{abstract}
\begin{thanks}{The authors were supported in part by the
Slovenian-Ukrainian research grant SLO-UKR 02-03/04.}
\end{thanks}

%\end{frontmatter}

\maketitle %\baselineskip15pt

One of the principal results linking Theory of Hyperspaces with
Theory of Retracts is\break \mbox{Wojdys\l awski} Theorem
\cite{Wo} asserting that the hyperspace $2^X$ of a compact metric
space $X$ is an absolute retract if and only if $X$ is a Peano
continuum (i.e., a continuous image of the interval $[0,1]$).
Among many characterizations of Peano continua let us recall the
Hahn-Mazurkiewicz-Sierpi\'nski Theorem (see \cite[\S50]{Ku})
asserting that a connected metrizable compact space is a Peano
continuum if and only if it is locally connected and the Bing
convexification theorem \cite{Bi}, \cite{Mo} characterizing Peano
continua as compacta admitting a convex metric. We recall that a
metric $d$ on $X$ is {\em convex} (resp. {\em almost convex}) if
for any $x,y\in X$ and positive reals $s,t$ with $d(x,y)\le s+t$
(resp. $d(x,y)<s+t$) there is $z\in X$ such that $d(x,z)\le s$ and
$d(z,y)\le t$. Each almost convex metric on a compact space is
convex. On the other hand, the standard metric on the space of
rational numbers is almost convex but fails to be convex.

Combining Bing's and Wojdyslawski's theorems we conclude that for
a compact space $X$ endowed with a convex metric $d$ the
hyperspace $2^X$ of non-empty closed subsets of $X$ is an absolute
retract. Generalizing this result to non-compact spaces,
C.Costantini and W.Kubi\'s \cite{CK} proved that for an almost
convex bounded metric space $X$ the hyperspace $2^X$ of all
non-empty closed subsets of $X$, endowed with the Hausdorff
metric, is an absolute neighborhood retract. In its turn
M.Kurihara, K.Sakai and M.Yaguchi \cite{KSY} showed that the
almost convexity of $X$ in the above result can be replaced by the
so-called uniform local $C^*$-connectedness (in the present paper
this property is called the uniform local chain
equi-connectedness).

In this paper we show that the latter result of \cite{KSY} can be
reversed. More precisely, the hyperspace $2^X$ of all non-empty
closed subsets of a metric space $X$ is an ANR (and AR) if and
only if $X$ is uniformly locally chain equi-connected (and chain
equi-connected). It is interesting to compare this result with
characterization of the ANR-property in other subspaces of $2^X$:
for a (connected) metric space $X$ the hyperspace $\mathcal F(X)$
of all non-empty finite subsets of $X$ is an ANR (an AR) if and
only if $X$ is locally path connected \cite{CN} while the
hyperspace $\mathcal K(X)$ of all non-empty compact subsets of $X$
is an ANR (an AR) if and only if $X$ is locally continuum
connected, see \cite{Cu}.
%
%In fact, the ANR-property in hyperspaces is equivalent to a wide
%spectrum of local properties having topological, uniform or metric
%nature. Before defining these local properties, recall some
%information related to hyperspaces.

Now it is time to give precise definitions. For a bounded metric
space $(X,d)$ let $2^X$ denote the hyperspace of all non-empty
closed subsets of $X$, endowed with the Hausdorff metric
$$d_H(A,B)=\max\{\sup_{x\in B}d(x,A),\sup_{x\in A}d(x,B)\}.$$ For
an unbounded metric space $(X,d)$ the Hausdorff metric can attain
infinite values but still determined a topology on $2^X$, called
the uniform topology. This topology depends not on a particular
metric $d$ on $X$ but on the uniformity generated by that metric,
see \cite[Ch.8]{En} for the theory of uniform spaces.

Each uniformity $\U$ generating the topology of a space $X$
induces the {\em Hausdorff uniformity} $2^\U$ on the hyperspace
$2^X$ of all non-empty closed subsets of $X$. This uniformity
$2^\U$ is generated by the base consisting of the entourages
$$2^U=\{(A,A')\in 2^X\times 2^X:A\subset B(A',U),\; A'\subset
B(A,U)\},\quad U\in\U.$$ Here, as expected, $B(A,U)=\bigcup_{a\in
A}B(a,U)$ where $B(a,U)=\{x\in X:(x,a)\in U\}$ is the $U$-ball
centered at $a\in X$, see \cite[8.5.16]{En}. The topology on $2^X$
induced by the Hausdorff uniformity $2^\U$ will be called the {\em
uniform topology} on $2^X$. Talking about topological properties
of $2^X$ we shall always refer to this (uniform) topology.

If the uniformity $\U$ on $X$ is generated by a metric $d$, then
the uniformity of the subspace $\Bd(X)\subset 2^X$ consisting of
all bounded non-empty closed subsets of $X$ is generated by the
Hausdorff metric $$d_H(A,B)=\max\{\sup_{x\in B}d(x,A),\sup_{x\in
A}d(x,B)\},$$ \cite[8.5.16.b]{En}. For arbitrary (non-necessarily
bounded) closed subsets $A,B$ of the metric space $(X,d)$ the
Hausdorff distance $d_H(A,B)$ can attain infinite values, but
letting $\rho_H=\min\{1,d_H\}$ we get a metric inducing the
uniformity $2^\U$ on the whole hyperspace $2^X$.

It should be mentioned that the uniform topology on $2^X$
coincides with the Vietoris topology if and only if $X$ is totally
bounded. On the other hand, these two topologies always coincide
on the subspace $\mathcal K(X)\subset 2^X$ of all non-empty
compact subsets of $X$, see \cite[8.5.16(c)]{En}.

The uniform as well as the Vietoris topologies on $2^X$ were
actively studied by general topologists, see \cite{En}, \cite{Be},
\cite{IN}. In spite of their coincidence on the hyperspace
$\mathcal K(X)$, for non-compact $X$ the uniform and Vietoris
topologies on $2^X$ differ substantially by their connectedness
properties. In particular, for a connected uniform space $X$ the
hyperspace $2^X$ endowed with the Vietoris topology is connected
but rarely is locally connected. On the other hand, the hyperspace
$2^X$ endowed with the uniform topology rarely is connected but
often is locally connected. A typical example of this phenomenon
is the hyperspace $2^\IR$ of the real line $\IR$. Endowed with the
Vietoris topology $2^\IR$ is connected but fails to be locally
connected. In contrast, endowed with the uniform topology the
hyperspace $2^\IR$ is locally connected but fails to be connected.

Striving to characterize metrizable uniform spaces $X$ whose
hyperspace $2^X$ is an ANR we invent that the ANR-property in
hyperspaces is equivalent to a wide spectrum of local properties
having topological, uniform, metric or extension nature.

%In this paper that the wide spectrum of local properties
%(beginning from being ANR and ending with the local chain
%connectedness) are equivalent for hyperspaces $2^X$ of metrizable
%uniform spaces. As expected, by a {\em metrizable uniform space}
%we understand a uniform space whose uniformity is generated by a
%metric (in this case the uniformity can be generated also by a
%bounded metric). At first we shall introduce several local
%connectedness properties of metric or uniform spaces. These
%properties fall into three categories having topological, uniform
%or metric nature.

We start with two local properties having topological nature. We
recall that a topological space $X$ is
\begin{itemize}
\item {\em locally path connected\/} [briefly (lpc)]
if for each $x_0\in X$ and a neighborhood $U\subset X$ of $x_0$
there is a neighborhood $V\subset X$ of $x_0$ such that each point
$x\in V$ can be linked with $x_0$ by a continuous path $f:[0,1]\to
U$ with $f(0)=x_0$, $f(1)=x$;
\item {\em locally connected\/} [briefly (lc)] if for each $x_0\in X$ and a
neighborhood $U\subset X$ of $x_0$ there is a connected subset
$C\subset U$, containing a neighborhood of $x_0$.
\end{itemize}

Next, we consider some properties having uniform nature. In some
cases defining certain property we shall simultaneously (in
parentheses) define its relative versions. For a point $x$ in a
metric space $(X,d)$ by $B(x,r)=\{y\in X:d(x,y)<r\}$ we denote the
open $r$-ball around $x$. We shall say that a metric space $(X,d)$
is
\begin{itemize}
\item {\em uniformly locally compact} if there is $\e>0$ such that each
closed $\e$-ball $\bar B(x,\e)=\{y\in X:d(x,y)\le\e\}$ is compact;
\item {\em chain connected\/} if for any points $x,y\in X$ and any
$\eta>0$ there is a sequence $x=x_0,x_1,\dots,x_l=y$ of points of
$X$ such that $d(x_i,x_{i-1})<\eta$ for all $i\le l$; such a
sequence $x=x_0,x_1,\dots,x_l=y$ is called an {\em $\eta$-chain}
linking $x$ and $y$ and $l$ is the {\em length} of this chain;
\item {\em chain equi-connected\/} if for any  $\eta>0$
there is a number $l\in\IN$ such that any points $x,y\in X$ can be
connected by an $\eta$-chain of length $\le l$;
\item {\em chain connected im kleinen\/} if there is  $\e>0$
such that any points $x,y\in X$ with $d(x,y)<\e$ can be linked by
an $\eta$-chain for any $\eta>0$;
\item {\em locally chain connected\/} [briefly (lcc)]
(at a subset $X_0\subset X$) if for each point $x_0$ in $X$ (in
$X_0$) and $\e>0$ there is $\delta>0$ such that for each $\eta>0$
and $x$ in $B(x_0,\delta)$ (in $X_0\cap B(x_0,\delta)$) there is
an $\eta$-chain $x_0,x_1,\dots,x_l=x$ of diameter $<\e$ linking
the points $x_0$ and $x$;
\item {\em uniformly locally  chain connected\/} [briefly (ulcc)]
if $\forall\e>0\;\exists\delta>0\;\forall\eta>0$ such that any
points $x,y\in X$ with $d(x,y)<\delta$ can be connected by an
$\eta$-chain of diameter $<\e$;
\item {\em uniformly locally chain equi-connected\/} [briefly
(ulcec)] (at a subset $X_0\subset X$) if $\forall
\e>0\;\exists\delta>0\;\forall\eta>0\;\exists l\in\IN$ such that
any points $x,y$ in $X$ (in $X_0$) with $d(x,y)<\delta$ can be
connected by an $\eta$-chain of diameter $<\e$ and length $\le l$;
\item {\em uniformly locally path connected\/} [briefly (ulpc)]
if for any $\e>0$ there is $\delta>0$ such that any points $x,y$
with $d(x,y)<\delta$ can be connected by a continuous path
$f:[0,1]\to X$ of diameter $<\e$ such that $f(0)=x$ and $f(1)=y$.
\item {\em locally trace connected\/} [briefly
(ltc)] if for each point $x_0\in X$ and $\e>0$ there is $\delta>0$
such that for any point $x\in X$ with $d(x,x_0)<\delta$ and a
countable dense subset $Q\supset\{0,1\}$ of $[0,1]$ there is a
uniformly continuous map $f:Q\to B(x_0,\e)$ such that $f(0)=x_0$
and $f(1)=x$; such a uniformly continuous map $f:Q\to X$ is called
a {\em trace} linking $x_0$ and $x$;
\item {\em uniformly locally trace connected\/} [briefly (ultc)]
if for any $\e>0$ there is $\delta>0$ such that for any countable
dense set $Q\supset\{0,1\}$ of $[0,1]$ any points $x,y$ with
$d(x,y)<\delta$ can be connected by a trace $f:Q\to X$ with $\diam
f(Q)<\e$.
\end{itemize}

Next, let us consider connectedness properties having metric
nature. By the {\em continuity modulus} of a function $f:X\to Y$
between metric spaces $(X,d)$, $(Y,\rho)$ we understand the
non-decreasing function $\w_f:(0,\infty)\to(0,\infty]$ defined by
$$\w_f(t)=\sup\{\rho(f(x),f(x')):x,x'\in X,\; d(x,x')\le t\}
\mbox{ for }t\in[0,\infty).$$ If the metric space $(X,d)$ is
almost convex, then the continuity modulus $\w_f$ is {\em
subadditive} in the sense that $\w_f(s+t)\le w_f(s)+\w_f(t)$. Let
us remark that the map $f$ is {\em uniformly continuous} if and
only if $\lim_{t\to+0}\w_f(t)=0$.

In the sequel by {\em a continuity modulus} we shall understand an
arbitrary non-decreasing function $\w:(0,\infty)\to(0,\infty)$
with $\lim_{t\to+0}\w(t)=0$.

We shall say that a metric space $(X,d)$ is
\begin{itemize}
\item {\em uniformly locally path equi-connected\/} [briefly
(ulpec)] if there are $\e_0>0$ and a continuity modulus $\w$ such
that any two points $x,y\in X$ with $d(x,y)<\e_0$ can be connected
by a path $f:[0,d(x,y)]\to X$ such that $f(0)=x$, $f(d(x,y))=y$,
and $\w_f\le \w$;
\item {\em uniformly locally trace equi-connected\/} [briefly
(ultec)] (at a subset $X_0\subset X$) if there are $\e_0>0$ and a
continuity modulus $\w$ such that for any two points $x,y$ in $X$
(in $X_0$) with $d(x,y)<\e_0$ and any countable dense subset
$Q\supset\{0,d(x,y)\}$ of $[0,d(x,y)]$ there is a uniformly
continuous function $f:Q\to X$ such that $f(0)=x$, $f(d(x,y))=y$,
and $\w_f\le \w$;
\item {\em $\w$-convex\/}, where $\w$ is a continuity modulus, if
for any points $x,y\in X$ and positive reals $s,t$ with
$d(x,y)<\w(s+t)$ there is a point $z\in X$ such that
$d(x,z)<\w(s)$ and $d(z,y)<\w(t)$.
\end{itemize}
Let us remark that a metric space $X$ is almost convex if and only
if it is $\w$-convex for a linear continuity modulus $\w(t)=at$.

By the {\em chain connected component} of a point $x_0$ in a
metric space $(X,d)$ we understand the set $C(x_0)=\{x\in X:$ for
any $\eta>0$ the points $x$ and $x_0$ can be linked by an
$\eta$-chain in $X\}$. We shall say that a uniform space $X$ is
\begin{itemize} \item {\em weakly convexifiable} [briefly (wcx)]
if $X$ is chain connected im kleinen and there is a concave
continuity modulus $\alpha$ such that for each concave continuity
modulus $\w\ge \alpha$ the uniformity of $X$ is generated by a
metric which is $\w$-convex on each chain connected component of
$X$.
\end{itemize}
We recall that a function $f:[0,\infty)\to[0,\infty)$ is {\em
concave} if $$f(ta+(1-t)b)\ge tf(a)+(1-t)f(b)$$ for any
$a,b\in[0,\infty)$ and $t\in[0,1]$.

Finally we define another three important local extension
properties of metric spaces which in case of hyperspaces are near
to the local connectedness properties discussed above. By a {\em
uniform neighborhood} of a set $X$ in a metric space $(M,d)$ we
understand any subset $U\subset M$ containing the
$\e$-neighborhood $B(X,\e)=\{x\in M:d(x,X)<\e\}$ of $X$ for some
$\e>0$. We shall say that a metric space $(X,d)$ is
\begin{itemize}
\item an {\em absolute neighborhood retract\/} [briefly ANR] if
for each metric space $M$ containing $X$ isometrically there is a
retraction $r:U\to X$ defined on some open neighborhood $U$ of $X$
in $M$;
\item an {\em absolute neighborhood uniform retract\/} [briefly ANUR]
if for each metric space $M$ containing $X$ isometrically there is
a uniformly continuous retraction $r:U\to X$ defined on some
uniform neighborhood $U$ of $X$ in $M$;
\item a {\em uniform absolute neighborhood retract\/} [briefly uANR]
if for each metric space $M$ containing $X$ isometrically there is
a retraction $r:U\to X$ defined on some uniform neighborhood $U$
of $X$ in $M$ and such that $r$ is uniformly continuous at $X$ in
the sense that for any $\e>0$ there is $\delta>0$ such that
$d(f(y),x)<\e$ for any $x\in X$ and $y\in U$ with $d(x,y)<\delta$.
\end{itemize}
Replacing in the above definition the neighborhood $U$ of $X$ by
the whole $M$ we will get the definitions of an {\em absolute
retract, absolute uniform retract}, and {\em uniform absolute
retract} [briefly, (AR), (AUR), and (uAR)].

It is interesting to note that the theories of ANUR's and uANR's
have different historical origins. ANUR's arose in Nonlinear
Functional Analysis \cite{BL} and were studies mainly by analysts
while uANR's first appeared in the topological paper of Michael
\cite{Mi} and were studied mainly by topologists, see \cite{Sa}.

We shall show that for each metric space $X$ the above-defined
properties relate as follows (the vertical arrows $\uparrow$
correspond to the implications holding under the additional
assumption of the {\em completeness} of $X$):

\begin{picture}(200,150)(-120,0)
\put(-1,130){(ANUR)}\put(46,130){$\Rightarrow$}
\put(60,130){(uANR)}\put(105,130){$\Rightarrow$}
\put(119,130){(ANR)}

\put(18,110){$\Downarrow$}\put(76,110){$\Downarrow$}\put(133,110){$\Downarrow$}

\put(3,90){(ulpec)}\put(44,90){$\Rightarrow$}
\put(65,90){(ulpc)}\put(105,90){$\Rightarrow$}
\put(124,90){(lpc)}
\put(13,70){$\Downarrow$}\put(22,70){$\uparrow$}
\put(73,70){$\Downarrow$}\put(82,70){$\uparrow$}
\put(130,70){$\Downarrow$}\put(139,70){$\uparrow$}
\put(150,82){\vector(1,-1){20}}

\put(-50,50){(wcx)}\put(-15,50){$\Leftrightarrow$}
\put(3,50){(ultec)}\put(44,50){$\Rightarrow$}
\put(65,50){(ultc)}\put(105,50){$\Rightarrow$}
\put(125,50){(ltc)}\put(175,50){(lc)}

\put(18,30){$\Updownarrow$}
\put(76,30){$\Downarrow$}\put(133,30){$\Downarrow$}
\put(170,44){\vector(-1,-1){20}}

\put(3,10){(ulcec)}\put(44,10){$\Rightarrow$}
\put(65,10){(ulcc)}\put(105,10){$\Rightarrow$}
\put(125,10){(lcc)}

\put(55,-15){Diagram 1}
\end{picture}
\bigskip
\bigskip

Except for the equivalences
(ulcec)$\Leftrightarrow$(ultec)$\Leftrightarrow$(wcx) (proven in
Section~\ref{sect2}) all the implications of Diagram 1 are rather
trivial (to see that the extension properties from the first row
of Diagram 1 imply the path-connectedness properties from the
second row, apply Eels-Arens-Wojdys\l awski Theorem (see
\cite[II.\S1]{BP}) asserting that each metric space is isometric
to a closed subset of a linear normed space).

It turns out that for the hyperspace $2^X$ of a metric space $X$
all these 14 properties (with possible exception of ANUR) are
equivalent.

Besides the hyperspace $2^X$ of a uniform space $X$ we shall be
interested in some its subspaces:
\begin{itemize}
\item $\Dis(X)=\{D\in 2^X: D$ is discrete$\}$
\item $\uDis(X)=\{D\in 2^X:D$ is uniformly discrete in $X\}$ and
\item $\Dis_0(X)=\{D\in \Dis(X):|D|\le\aleph_0$ and
$D$ is totally bounded or uniformly
discrete in $X\}$.
\end{itemize}
We call a subset $D$ of a metric space $(X,d)$ {\em uniformly
discrete} if there is $\e>0$ such that $|F\cap B(x,\e)|\le 1$ for
each point $x\in X$. It is easy to see that a subset $F\subset X$
is uniformly discrete in $X$ if and only if $F$ is {\em
$\e$-separated} for some $\e>0$ in the sense that $d(x,y)\ge\e$
for all distinct $x,y\in F$. Using the Zorn Lemma it is easy to
show that the set $\uDis(X)$ is dense in $2^X$. Note also that for
a complete metric space $X$ the set $\Dis_0(X)$ coincides with the
collection of all at most countable uniformly discrete subsets of
$X$ (because closed discrete totally bounded subsets in complete
uniform spaces are finite).

The main result of this paper is the following characterizing

\begin{theorem}\label{t1} For a metric space $(X,d)$ and its hyperspace
$2^X$ endowed with the ``Hausdorff'' metric $\min\{1,d_H\}$ the
following conditions are equivalent:
\begin{enumerate}
\item $2^X$ is an ANR;
\item $2^X$ is locally connected;\smallskip

\item $2^X$ is a uniform ANR;
\item $2^X$ is uniformly locally path equi-connected;
\item $2^X$ is locally chain connected;
\item $2^X$ is locally chain connected at $\Dis_0(X)$;
\smallskip

\item $X$ is uniformly locally chain equi-connected;
\item $X$ is uniformly locally trace equi-connected;
\item $X$ is weakly convexifiable.
\medskip

\noindent \hskip-30pt Moreover, if the space $X$ is complete, then
the conditions \textup{(1)--(9)} are equivalent to
\smallskip

\item $X$ is uniformly locally path equi-connected;
\medskip

\noindent \hskip-30pt If $X$ is uniformly locally compact, then
the conditions \textup{(1)--(10)} are equivalent to
\smallskip

\item Each space $\Hip\subset 2^X$ containing $\Dis(X)$ is a uniform ANR;
\smallskip

\noindent \hskip-30pt If $X$ is chain equi-connected, then the
conditions \textup{(1)--(9)} are equivalent to
\smallskip

\item $2^X$ is an absolute retract;
\item $2^X$ is a uniform absolute retracts.
\end{enumerate}
\end{theorem}

The implication $(7)\Rightarrow(3)$ was first proved in \cite{KSY}
while some particular cases of $(10)\Rightarrow(11)$ were proved
by M.Kurihara, see \cite{KSY}. By Proposition~4.6 of \cite{CK},
the hyperspace $2^X$ of a metric space $X$ is chain connected if
and only if $X$ is chain equi-connected. Combining this result
with Theorem~\ref{t1} we get characterizing

\begin{theorem}\label{t1c} The hyperspace $2^X$ of a metric space $X$ is an
absolute retract if and only if $2^X$ is a uniform absolute
retract if and only if $X$ is chain equi-connected and uniformly
locally chain equi-connected.
\end{theorem}

Some implications of Theorem~\ref{t1} hold in a more general
setting. For a closed subset $F$ of a metric space $X$ let
$\Dis_0(F)=\{D\in\Dis_0(X):D\subset F\}$.

\begin{theorem}\label{t2} For a metric space $X$ and
an open subspace $\Hip\subset 2^X$ such that $\Dis_0(F)\subset
\mathcal H$ for each $F\in\mathcal H$ the following conditions are
equivalent:
\begin{enumerate}
\item $\Hip$ is an ANR;
\item $\Hip$ is locally path connected;
\item $\Hip$ is locally chain connected;
\item $2^X$ is locally chain connected at $\Dis_0(F)$ for each $F\in\mathcal H$;
\item $X$ is uniformly locally chain equi-connected at each $F\in\Hip$;
\item $X$ is uniformly locally trace equi-connected at each $F\in\Hip$;
\end{enumerate}
\end{theorem}

Applying Theorem~\ref{t2} to the hyperspace $\mathcal
H=\Bd(X)\subset 2^X$ of all closed bounded subsets of a metric
space $X$ we get

\begin{corollary}\label{t3} For the hyperspace $\Bd(X)\subset 2^X$ of all
non-empty closed bounded subsets of a metric space $X$ the
following conditions are equivalent:
\begin{enumerate}
\item $\Bd(X)$ is an ANR;
\item $\Bd(X)$ is locally path connected;
\item $\Bd(X)$ is locally chain connected;
\item $\Bd(X)$ is locally chain connected at $\Bd(X)\cap\Dis_0(X)$;
\item $X$ is uniformly locally chain equi-connected
at each bounded subset of $X$;
\item $X$ is uniformly locally trace equi-connected
at each bounded subset of $X$;
\end{enumerate}
\end{corollary}

Theorem~\ref{t3} allows us to characterize the AR-property in the
hyperspaces $\Bd(X)$.

\begin{theorem}\label{t4} The hyperspace $\Bd(X)$ of a metric
space $X$ is an absolute retract if and only if
\begin{enumerate}
\item[(i)] each bounded subset of $X$ lies in a bounded chain equi-connected subspace of $X$ and
\item[(ii)] $X$ is uniformly locally chain
equi-connected at each bounded subset of $X$.
\end{enumerate}
\end{theorem}

\begin{problem} Characterize metric spaces whose hyperspace
is an absolute (neighborhood) {\em uniform} retract.
\end{problem}

The above characterizations imply several unexpected corollaries.

\begin{corollary}\label{cor2} Let $X$ be a dense subset of a metric space $M$.
\begin{enumerate}
\item The hyperspace $2^X$ is an absolute (neighborhood) retract
if and only if so is the hyperspace $2^M$.
\item The hyperspace $\Bd(X)$ is an absolute (neighborhood) retract
if and only if so is the hyperspace $\Bd(M)$.
\end{enumerate}
\end{corollary}

Another unexpected corollary is an amusing characterization of
normable spaces. We recall that a linear topological space $X$ is
{\em normable} if its topology is determined by a norm.

\begin{corollary}\label{cor3} A metrizable locally convex space $X$
is normable if and only if its hyperspace $2^X$ is an ANR.
\end{corollary}

\begin{proof} Each normed space $X$ carries a convex metric and
consequently is uniformly locally chain equi-connected. Applying
Theorem~\ref{t1}, we conclude that the hyperspace $2^X$ is an ANR.

Assume conversely that the hyperspace $2^X$ of a metrizable
locally convex space $X$ is an ANR. Then Theorem~\ref{t1} implies
that $X$ is uniformly locally chain equi-connected. This means
that there is a convex neighborhood $U$ of the origin of $X$ such
that for any convex neighborhood $W\subset X$ of the origin there
is $l\in\IN$ such that for any $x\in U$ there is a chain
$0=x_0,\dots,x_l=x$ with $x_i-x_{i-1}\in W$ for all $i\le l$. This
implies that $\underset{l}{\underbrace{W+\dots+W}}=l\,W\supset U$,
which means that the set $U$ is bounded in $X$. Then $X$ contains
a bounded convex neighborhood $U$ of the origin and hence $X$ is
normable, see \cite[II.2.1]{Sch}.
\end{proof}

This corollary implies that the hyperspace $2^{\IR^\w}$ of
$\IR^\w$, the countable product of lines, fails to be an ANR (in
spite of the fact that $\IR^\w$ is an absolute uniform retract).
The first example of an metric absolute retract $X$ whose
hyperspace $2^X$ fails to be an ANR was constructed in \cite{KSY}.

\begin{problem} Characterize metric linear spaces whose
hyperspaces are ANR's.
\end{problem}

Finally let us pose an intriguing open problem related to the Bing
convexification Theorem \cite{Bi} and the implication
(2)$\Rightarrow$(9) of Theorem~\ref{t1}. According to this
implication, the uniformity of any chain connected metric space
$X$ with locally connected hyperspace $2^X$ is generated by a
$\w$-convex metric for some concave continuity modulus $\w$. Let
us remark that a metric $d$ is almost convex if and only if it is
$\w$-convex for a linear continuity modulus $\w(t)=at$.

\begin{problem} Assume that the hyperspace $2^X$ of some chain
connected metric space $X$
is an ANR. Is the uniformity of $X$ generated by an almost convex
metric?
\end{problem}

The answer to this problem is affirmative if $X$ is totally
bounded.

\begin{theorem} The hyperspace $2^X$ of a totally bounded
metric space is an absolute retract if and only if the uniformity
of $X$ is generated by an almost convex metric.
\end{theorem}

\begin{proof} The ``if'' part of this theorem was proved in \cite{CK}
(and can be derived from Theorem~\ref{t1c}). To prove the ``only
if'' part, assume that $2^X$ is an AR. Since $X$ is totally
bounded, the completion $Y$ of $X$ is compact. Since $X$ is dense
in $Y$, Corollary~\ref{cor2} implies that the hyperspace $2^Y$ is
an AR. Combining Wojdys\l awski's and Bing's Theorems we conclude
that the uniformity of $Y$ is generated by a convex metric $d$.
Then $d$ restricted to $X$ is almost convex and generates the
uniformity of $X$.
\end{proof}
\medskip

Now let us pass to the proofs of our results.

\section{Proofs of the implications from Diagram 1}\label{sect2}

In this section we shall prove the non-trivial equivalences
(ulcec)$\Leftrightarrow$(ultec)$\Leftrightarrow$(wcx) from Diagram
1.

\begin{lemma}\label{lem1} If a metric space $(X,d)$
is uniformly locally trace equi-connected at a subset $X_0\subset
X$, then $X$ is uniformly locally chain equi-connected at $X_0$.
\end{lemma}

\begin{proof} Assuming
that $X$ is uniformly locally trace equi-connected at $X_0$, find
$\e_0>0$ and a continuity modulus $\w$ such that any two points
$x,y\in X_0$ with $d(x,y)<\e_0$ can be linked by a trace $f:Q\to
X$ defined on a dense subset $Q\supset\{0,d(x,y)\}$ of
$[0,d(x,y)]$ and such that $f(0)=x$, $f(d(x,y))=y$, and $\w_f\le
\w$.

To prove that $X$ is uniformly locally chain equi-connected in
$X_0$, fix any $\e>0$ and find $0<\delta<\e_0$ such that
$\w(\delta)<\e$. Next, given $\eta>0$ find $l\in\IN$ such that
$\w(\delta/l)<\eta$. Fix any points $x,y\in X_0$ with
$d(x,y)<\delta$. It follows that there is a uniformly continuous
function $f:Q\to X$ defined on a dense subset
$Q\supset\{0,d(x,y)\}$ of $[0,d(x,y)]$ such that $f(0)=x$,
$f(d(x,y))=y$, and $\w_f\le \w$. Since $d(x,y)<\delta$, we can
find points $0=t_0\le t_1\le\dots\le t_l=d(x,y)$ in $Q$ such that
$t_i-t_{i-1}<\delta/l$. Letting $x_i=f(t_i)$ for $i\le l$ we will
get an $\eta$-chain $x_0,\dots,x_l$ of length $\le l$ linking the
points $x=x_0$ and $y=x_l$. Since $|t_i-t_j|<\delta$, we get
$d(x_i,x_j)=d(f(t_i),f(t_j))\le\w_f(|t_i-t_j|)\le\w(\delta)<\e$,
which means that this chain has diameter $<\e$.
\end{proof}

Lemma~\ref{lem1} can be partly reversed.

\begin{lemma}\label{lem2} A metric space $(X,d)$ is
uniformly locally trace equi-connected at a subset $X_0\subset X$,
provided $X$ is uniformly locally chain equi-connected at some
$r$-neighborhood $B(X_0,r)$ of $X_0$.
\end{lemma}

\begin{proof} Assume that $X$ is uniformly locally chain
equi-connected at $B(X_0,r_0)$ for some $r_0>0$. To prove the
lemma it suffices to find $\e_0>0$ and a continuity modulus $\w$
such that for any points $x,y\in X_0$ with $d(x,y)<\e_0$ there is
a uniformly continuous function $f:Q\to X$ defined on {\em some}
dense subset $Q\supset\{0,d(x,y)\}$ of $[0,d(x,y)]$ and such that
$\w_f\le\w$, $f(0)=x$ and $f(d(x,y))=y$. Having such a function
$f:Q\to X$ and given any countable dense subset
$Q'\supset\{0,d(x,y)\}$ of $[0,d(x,y)]$ we can find an increasing
Lipschitz homeomorphism $h$ of $[0,d(x,y)]$ with Lipschitz
constant 2 such that $h(Q')=Q$. Then the composition $f\circ
h|Q':Q'\to X$ is a uniformly continuous map such that $f\circ
h(0)=x$, $f\circ h(d(x,y))=y$ and $\w_{f\circ h}(t)\le \w(2t)$ for
all $t\ge 0$.

Using the definition of (ulcec) of $X$ at $X_0$ we can construct a
decreasing sequence $(\delta_n)_{n\in\IN}$ of positive reals such
that for each $n\in\IN$ and $\eta>0$ there is a number
$l=l(n,\eta)\in\IN$ such that any two points $x,y\in U$ with
$d(x,y)<\delta_n$ can be linked by an $\eta$-chain of diameter
$<2^{-n}r_0$ and length $\le l$. Without loss of generality,
$\delta_n\le 2^{-n}r_0$ for all $n$.

Now for every $n\in\IN$ take any $l_n\ge l(n,\delta_{n+1})$ and
let $\Delta_n=1/(l_1\cdots l_n)$. Let also $\Delta_0=1$. Replacing
$(l_n)$ by a larger sequence, if necessary, we can assume that
$\Delta_{n}<\delta_{n+2}$ for all $n\in\IN$.

For every $m\in\IN$ consider the finite subset
$$Q_m=\big\{\sum_{k=1}^m\frac{i_k}{\Delta_k}:0\le i_k<l_k \mbox{
for all $k\le m$}\}$$ of $[0,1)$ and let $Q=\bigcup_{m\ge1}Q_m$.
It is clear that $Q$ is a countable dense subset of $[0,1]$.

A continuity modulus $\w$ will be defined as the supremum
$\w=\sup_{n\in\IN}\w_n$ of the sequence $(\w_n)$ of continuity
moduli defined recursively: $\w_0\equiv 0$ and
$$\w_{n+1}(t)=\begin{cases} 0&\mbox{ if $t<\Delta_{n+1}$};\\
2^{-n}r_0+\w_n(t)&\mbox{ if $t\ge\Delta_{n+1}$}
\end{cases}
$$ for $n\in\IN$. It easy to see that each function $\w_n$ is
non-decreasing and the supremum $\w=\sup_{n\in\IN}\w_n$ is a
well-defined continuity modulus.

To finish the proof it rests to verify that $\e_0=\delta_2$ and
$\w$ satisfy our requirements. Take any two points $x,y\in X_0$
with $d(x,y)<\e_0=\delta_2$. Find $k\ge 2$ with $\delta_{k+1}\le
d(x,y)<\delta_k$. Then $\Delta_{k-1}<\delta_{k+1}\le d(x,y)$.
%Let $I=[0,\Delta_k]$.

By induction we shall construct a sequence of functions
$\{f_m:Q_m\to B(X_0,r_0)\}_{m\ge k-1}$ such that the following
conditions are satisfied for every $m\ge k-1$:
\smallskip

1) $f_m(q)=y$ for all $q\in Q_m\cap[\Delta_{k-1},1)$;

2) $f_{m+1}|Q_{m}=f_m$;

3) $f_m(Q_m)\subset B(X_0,r_0\sum_{i=k}^{m}2^{-i})$;

4) $d(f_m(p),f_m(q))<\delta_{m+1}$ for any neighbor points $p,q\in
Q_m$;

5) $\diam f_m(Q_m\cap[p,q])< 2^{-m}r_0$ for any neighbor points
$p,q$ of $Q_{m-1}$.

6) $\w_{f_m}\le \w_m$.
\medskip

Two points $p,q$ of a finite subset $F\subset [0,1]$ are called
{\em neighbor points} if $(p,q)$ is a connected component of
$[0,1]\setminus F$.

To start the inductive construction let $f_{k-1}(0)=x$ and
$f_{k-1}(q)=y$ for any $q\in Q_{k-1}\cap [\Delta_{k-1},1)$. Assume
that for some $m\ge k-1$ functions $f_{k-1},\dots,f_m$ satisfying
the conditions (1)--(6) have been constructed. Take any neighbor
points $p,q$ of $Q_m$. If $\min\{p,q\}\ge \Delta_{k-1}$ let
$f_{m+1}|[p,q]\cap Q_{m+1}=\{y\}$. Otherwise use the conditions
(3) and (4) to conclude that the points $f_m(p),f_m(q)$ belong to
$B(X_0,r_0)$ and satisfy $d(f_m(p),f_m(q))<\delta_{m+1}$. The
choice of the sequences $(\delta_i)$ and $(l_i)$ guarantees the
existence of a $\delta_{m+2}$-chain
$f_m(p)=c_0,\dots,c_{l_{m+1}}=f_m(q)$ of diameter $<r_0/2^{m+1}$
and length $l_{m+1}$ connecting the points $f_m(p)$ and $f_m(q)$.
The intersection $(p,q)\cap Q_{m+1}$ contains exactly $m_{m+1}-1$
points, so we can assign to each of these points a point $c_i$
from the chain to satisfy the condition (4). Since $\diam
f_{m+1}(Q_{m+1}\cap[p,q])<2^{-(m+1)}r_0$, we see that the
conditions (3,5) are satisfied.

To show that $\w_{f_{m+1}}\le\w_{m+1}$ it suffices to verify that,
$d(f_{m+1}(p),f_{m+1}(q))\le \w_{m+1}(|p-q|)$ for any two distinct
points $p,q\in Q_{m+1}$. Given such points $p,q\in Q_{m+1}$ find
points $p',q'\in Q_m$ with $\max\{|p-p'|,|q-q'|\}<\Delta_m$ and
$|p'-q'|\le|p-q|$. Then $$
\begin{aligned} d(f_{m+1}(p),f_{m+1}(q))&\le
d(f_{m+1}(p),f_m(p'))+d(f_{m+1}(q),f_m(q'))+d(f_m(p'),f_m(q'))\le\\
&\le2\cdot 2^{-(m+1)}r_0+\w_{m}(|p'-q'|)\le
2^{-m}r_0+\w_m(|p-q|)=\w_{m+1}(|p-q|).
\end{aligned}
$$

This completes the inductive step. Now let
$f=\bigcup_{m\in\IN}f_m:Q\to B(X_0,r_0)$ and observe that
$f(0)=x$, $f(\Delta_{k-1})=\lim_{q\to d(x,y)}f(q)=y$ and
$\w_f\le\sup\w_m=\w$.
\end{proof}

\begin{lemma}\label{lem3} A metric space $(X,d)$
is uniformly locally chain equi-connected, provided $X$ is chain
connected im kleinen and the uniformity of $X$ is generated by a
metric $d$ which is $\w$-convex on each chain connected component
of $X$ for some continuity modulus $\w$.
\end{lemma}

\begin{proof} To show that
$X$ is uniformly locally chain equi-connected, fix arbitrary
$\e>0$. Replacing $\e$ by a smaller positive number we can assume
that any two points $x,y\in X$ with $d(x,y)<\e$ belong to the same
chain connected component of $X$.

For this $\e>0$ find $a>0$ such that $\w(a)<\frac{\e}2$ and let
$\delta=\w(a)$. Given $\eta>0$ find $l\in\IN$ such that
$\w(a/l)<\eta$.

Now take any points $x,y\in X$ with
$d(x,y)<\delta=\w(a)<\frac\e2$. They belong to some chain
connected component $C$ on which the metric $d$ is $\w$-convex.
Let $x_0=x$. Since $d$ is $\w$-convex on $C$ and
$d(x,y)<\delta=\w(a)=\w(\frac al+ \frac{a(l-1)}l)$ there is a
point $x_1\in C$ such that $d(x_0,x_1)<\w(a/l)<\eta$ and
$d(x_1,y)<\w(\frac{a(l-1)}l)$. Applying the $\w$-convexity of $d$
on $C$ once more, we shall find a point $x_2\in C$ such that
$d(x_1,x_2)<\w(a/l)<\eta$ and $d(x_2,y)<\w(\frac{a(l-2)}l)$.
Proceeding in this way we shall construct a sequence
$x=x_0,x_1,\dots,x_l=y$ of points of $C$ such that
$d(x_i,x_{i+1})<\w(a/l)<\eta$ and $d(x_i,y)<\frac{a(l-i)}l$ for
all $i<l$. It is clear that $x_0,\dots,x_l$ is an $\eta$-chain of
length $l$ and diameter $<2\w(a)<\e$ linking the points $x$ and
$y$.
\end{proof}

In the proof of convexification Lemma~\ref{lem5} we will exploit
the following elementary fact which can be proven by a simple
geometric argument.

\begin{lemma}\label{lem4} If $f$ is a concave continuity modulus, then
\begin{enumerate}
\item $f(b+d)-f(b)\le
f(a+d)-f(a)$ for any $0\le a\le b$ and $d\ge 0$;
\item $f(s+t)-f(s)-f(t)\le f(a+b)-f(a)-f(b)$ for any $0\le a\le s$ and
$0\le b\le t$.
\end{enumerate}
\end{lemma}

\begin{lemma}\label{lem5} If a metric space $(X,d)$ is uniformly locally chain
equi-connected, then $X$ is weakly convexifiable.
\end{lemma}

\begin{proof} If $X$ is (ulcec), then it is chain connected im kleinen and
is (ultec) by Lemma~\ref{lem2}. Consequently there are $R>0$ and a
continuity modulus $\w$ such that any points $x,y\in X$ with
$d(x,y)\le R$ can be linked by a uniformly continuous function
$f:Q\to X$ defined on a dense subset $Q\supset\{0,d(x,y)\}$ of
$[0,d(x,y)]$ and such that $f(0)=x$, $f(d(x,y))=y$ and $\w_f\le
\w$.

Let $C$ be the closed convex hull of the set
$\{(t,\w(t)):t\in(0,R]\}\cup\{(nR,(n+2)\w(R)):n\in\IN\}$ in the
real plane $\IR^2$. It can be easily shown that the function
$\gamma(t)=\max\{y\in C: (t,y)\in C\}$ is a concave continuity
modulus with $\gamma(t)\ge \w(t)$ for all $t\le R$ and
$\gamma(iR)\ge (i+2)\w(R)$ for all $i\in\IN$. Then $2\gamma$ is a
concave continuity modulus as well.

Given a concave continuity modulus $\alpha\ge 2\gamma$ we shall
define a metric $\rho$ on $X$ as follows. For a pair of points
$(x,y)\in X$ let $R(x,y)$ be the set of all positive real numbers
$r$ for which there is a uniformly continuous map $f:Q\to X$ from
a countable dense subset $Q\supset\{0,r\}$ of $[0,r]$ such that
$f(0)=x$, $f(r)=y$ and $\w_f\le \alpha$. It is easy to see that
$R(x,y)=R(y,x)$ for all $x,y\in X$.

Let us show that $R(x,y)\ne\emptyset$ for any points $x,y$ from
the same chain connected component. Given such points $x,y$ find
an $R$-chain $x=x_0,\dots,x_l=y$ linking points $x,y$. For every
$i\le l$ we get $d(x_{i-1},x_i)<R$. Consequently, we can find a
uniformly continuous function $f_i:Q_i\to X$ defined on a
countable dense subset $Q_i\supset\{(i-1)R,iR\}$ of the interval
$[(i-1)R,iR]$ such that $f((i-1)R)=x_{i-1}$, $f(iR)=x_i$ and
$w_{f_i}\le \w$. Let $Q=\bigcup_{i\le l}Q_i$ and $f=\bigcup_{i\le
l}f_i:Q\to X$. It is clear that $f$ is a uniformly continuous
function with $f(0)=x$ and $f(lR)=y$.

Let us show that $w_f(t)\le \alpha(t)$ for all $t\ge0$. Take any
points $q<q'$ in $Q$ with $|q-q'|\le t$. Fix numbers $i\le j\le l$
such that $(i-1)R<q\le iR$ and $(j-1)R\le q'< jR$. If $i=j$, then
$d(f(q),f(q'))\le \w_{f_i}(t)\le\w(t)\le\gamma(t)\le\alpha(t)$. If
$j=i+1$, then $$\begin{aligned} d(f(q)),f(q'))&\le
d(f_i(q),f_i(iR))+d(f_j(iR),f_j(q'))\le\\&\le
\w_{f_i}(iR-q)+\w_{f_j}(q'-iR)\le2\w(t)\le2\gamma(t)\le\alpha(t).
\end{aligned}$$
If $j>i+1$, then $$ (j-i-1)R\le q'-q=(q'-(j-1)R)+(j-i-1)R+
(iR-q)<2R+(j-i-1)R=(j-i+1)R$$ and $$\begin{aligned}
d(f(q),f(q'))&\le
d\big(f(q),f(iR)\big)+d\big(f(iR),f((i+1)R)\big)+\dots+d\big(f((j-1)R),f(q')\big)\le\\&\le
(j-i+1)\w(R)\le\gamma((j-i-1)R)\le
\gamma(q'-q)\le\gamma(t)\le\alpha(t).\end{aligned}$$

For points $x,y\in X$ let $r(x,y)=\inf (\{\infty\}\cup R(x,y))$.
Observe that $r(x,y)=r(y,x)$ and $r(x,y)<\infty$ if and only if
the points $x,y$ belong to the same chain connected component of
$X$. It follows from the definition of $r$ that for any $x,y\in X$
we get $d(x,y)\le \alpha\circ r(x,y)$. Moreover, $r(x,y)\le
d(x,y)$ if $d(x,y)\le R$.

Define a function $p:X\to[0,\infty]$ letting
$$p(x,y)=\inf\big\{\sum_{i=1}^m \alpha\circ r(x_{i-1},x_i):
x=x_0,\dots,x_m=y\big\}.$$ It is clear that this function is
symmetric and satisfies the triangle inequality (we assume that
$x+\infty=\infty\ge x$ for any $x\in[0,\infty]$). Also
$p(x,y)\le\alpha\circ r(x,y)\le \alpha\circ d(x,y)$ if $d(x,y)\le
R$.

On the other hand, for any sequence $x=x_0,\dots,x_m=y$ in $X$ we
get $$\sum_{i=1}^m\alpha\circ r(x_{i-1},x_i)\ge \sum_{i=1}^m
d(x_{i-1},x_i)\ge d(x_0,x_m)=d(x,y)$$ which implies that
$d(x,y)\le p(x,y)$, and $p(x,y)\le\alpha\circ d(x,y)$ if
$d(x,y)\le R$.

Take any subset $S\subset X$ meeting each chain connected
component of $X$ in a unique point. Define a metric $\rho$ on $X$
letting $\rho(x,y)=p(x,y)$ if $x,y$ belong to the same connected
component and $\rho(x,y)=2R+p(x,x_0)+p(y,y_0)$ where $x_0\in S$
(resp. $y_0\in S$) belongs to the chain connected component of $x$
(resp. $y$). It is easy to see that $\rho$ is a metric. Moreover,
$d(x,y)\le\rho(x,y)=p(x,y)\le \alpha\circ d(x,y)$ if $d(x,y)<R$.
This means that the metric $\rho$ is uniformly equivalent to $d$
and thus generates the uniformity of $X$.

It rests to verify that $\rho$ is $\alpha$-convex on each
connected component $C$ of $X$. Fix any points $x,y\in C$ and pick
positive real numbers $s,t$ with $\rho(x,y)<\alpha(s+t)$. If
$\rho(x,y)<\alpha(s)$, then put $z=y$ and note that
$\rho(x,z)=\rho(x,y)<\w(s)$ while $\rho(z,y)=0<\w(t)$. So further
we assume that $\rho(x,y)\ge w(s)$.

By the definition of $p(x,y)=\rho(x,y)$ there is a chain
$x=x_0,x_1,\dots,x_m=y$ such that
$\rho(x,y)\le\sum_{i=1}^m\alpha\circ r(x_{i-1},x_i)<\alpha(s+t)$.
Let $j\ge0$ be the largest number such that
$\sum_{i=1}^j\alpha\circ r(x_{i-1},x_i)<\alpha(s)$. The inequality
$\rho(x,y)\ge \alpha(s)$ implies that $j<l$.

Let $A=\sum_{i=1}^j\alpha\circ r(x_{i-1},x_i)$,
$B=\sum_{i=j+2}^l\alpha\circ r(x_{i-1},x_i)$, and fix any $\Delta>
r(x_{j-1},x_j)$ with
\begin{equation}\label{eq1}
A+\alpha(\Delta)+B<\alpha(s+t).
\end{equation}
It follows from the definition of $r(x_{j-1},x_j)<\Delta$ that
there is a uniformly continuous function $f:Q\to X$ defined on a
countable dense subset $Q\supset\{0,\Delta\}$ of $[0,\Delta]$ and
such that $f(0)=x_{j-1}$, $f(\Delta)=x_j$ and $w_f\le\alpha$.

By the choice of $j$, we get $A<\alpha(s)$, and
$A+\alpha(\Delta)\ge \alpha(s)$. Using the continuity of $\alpha$
find $\delta\in[0,\Delta]$ such that $A+\alpha(\delta
)=\alpha(s)$. It follows that $\delta \le s$ and
\begin{equation}\label{eq2}
B<\alpha(s+t)-\alpha(\Delta)-A=
\alpha(s+t)-\alpha(\Delta)-(\alpha(s)-\alpha(\delta )).
\end{equation}
We claim that $\Delta-\delta <t$. Assuming the converse we would
get
\begin{equation}\label{eq3}\alpha(\Delta)-\alpha(\delta )=\alpha(\delta +(\Delta-\delta ))-\alpha(\delta )\ge
\alpha(\delta +t)-\alpha(\delta )\ge
\alpha(s+t)-\alpha(s).\end{equation} The last inequality follows
from $\delta \le s$ and Lemma~\ref{lem4}(1). Combining
(\ref{eq1}--\ref{eq3}) we will get
$\alpha(s+t)>A+\alpha(\Delta)\ge (\alpha(s)-\alpha(\delta
))+\alpha(\Delta) \ge(\alpha(s)-\alpha(\delta
))+\alpha(s+t)-\alpha(s)+\alpha(\delta )=\alpha(s+t), $ which is a
contradiction. Thus $\delta \le s$ and $\Delta-\delta < t$ and we
can apply Lemma~\ref{lem4}(2) to conclude that $\alpha(\delta
)+\alpha(\Delta-\delta
)-\alpha(\Delta)\le\alpha(s)+\alpha(t)-\alpha(s+t)$ and hence
$$\alpha(\Delta-\delta )\le
(\alpha(s)+\alpha(t)-\alpha(s+t))+(\alpha(\Delta)-\alpha(\delta
)).$$ Adding to this inequality the inequality (\ref{eq2}), we
will get $$\begin{aligned}\alpha(\Delta-\delta )+B&<
\alpha(\Delta-\delta )+(\alpha(s+t)-\alpha(s))-
(\alpha(\Delta)-\alpha(\delta ))\le\\ &\le
\big(\alpha(s)+\alpha(t)-\alpha(s+t))+(\alpha(\Delta)-\alpha(\delta
)\big)+\\&+
\big(\alpha(s+t)-\alpha(s)\big)-\big(\alpha(\Delta)-\alpha(\delta
)\big)=\alpha(t).\end{aligned}$$ By the continuity of $\alpha$
find $q\in Q\cap[0,\delta )$ such that
$\alpha(\Delta-q)+B<\alpha(t)$. Then $A+\alpha(q)<A+\alpha(\delta
)=\alpha(s)$ and $B+\alpha(\Delta-q)<\alpha(t)$. It can be shown
that for the point $z=f(q)$ we get $\rho(x,z)=p(x,z)\le
A+\alpha(q)<\alpha(s)$ and $\rho(z,y)=p(y,z)\le
B+\alpha(\Delta-q)<\alpha(t)$. This completes the proof of the
$\alpha$-convexity of the metric $\rho$ on the chain connected
component $C$.
\end{proof}

\section{Lawson semilattices and (uniform) ANR's}

By a {\em topological semilattice} we understand a pair $(X,\vee)$
consisting of a topological space $X$ and a continuous associative
commutative idempotent operation $\vee:X\times X\to X$. A
topological semilattice $(X,\vee)$ is called a {\em Lawson
semilattice} if $X$ has a base of the topology consisting of
subsemilattices. The following observation made in \cite{BKS}
allows us to reduce the study of the ANR-property in Lawson
semilattices to verifying the local path connectedness.

\begin{proposition}\label{prop1} A metrizable Lawson semilattice
$L$ is an ANR (an AR) if and only if $L$ is locally path-connected
(and connected).
\end{proposition}

This result has a uniform counterpart. By a {\em Lawson metric
semilattice} we shall understand a metric space $(X,d)$ endowed
with a semilattice operation $\vee:X\times X\to X$ such that\break
$d(x\vee y,x'\vee y')\le \max\{d(x,x'),d(y,y')\}$ for all
$x,y,x',y'\in X$. It is easy to see that each Lawson metric
semilattice is Lawson as a topological semilattice. The following
important result proven in \cite{KSY} is a uniform counterpart of
Proposition~\ref{prop1}.

\begin{proposition}\label{prop2}
A Lawson metric semilattice $X$ is a uANR (and uAR) if and only if
it is uniformly locally path-connected (and chain connected).
\end{proposition}

It is easy to see that for any bounded metric space $(X,d)$ the
hyperspace $(2^X,d_H,\cup)$ is a Lawson metric semilattice with
respect to the operation of union $\cup$. Hence $2^X$ is a uANR
(and a uAR) if and only if it is uniformly locally path-connected
(and chain connected).

Propositions~\ref{prop1}, \ref{prop2} show that for Lawson metric
semilattices the implications (ANR)$\Rightarrow$(lpc) and
(uANR)$\Rightarrow$(ulpc) can be reversed. We do not known if this
can be done for the implication (ANUR)$\Rightarrow$(ulpec).

\begin{problem} Suppose that a Lawson metric semilattice $X$
is uniformly locally path equi-connected. Is $X$ an absolute
neighborhood uniform retract? Is $2^X$ an ANUR?
\end{problem}

\section{ANR-property of hyperspaces $\Dis(X)$}

The aim of this section is to prove the implication
$(10)\Rightarrow(11)$ of Theorem~\ref{t1}. We shall say that a
metric space $(X,d)$ has {\em totally bounded small balls} if
there is $r>0$ such that each ball $B(x,r)$, $x\in X$, is totally
bounded. It is clear that each uniformly locally compact space has
totally bounded small balls.

\begin{proposition}\label{prop3} Let $(X,d)$ be a uniformly locally path
equi-connected metric space. If $X$ has totally bounded small
balls, then the hyperspace $\Dis(X)$ is a uniform ANR.
\end{proposition}

In light of Proposition~\ref{prop2} it suffices to prove that
$\Dis(X)$ is uniformly locally path connected. The latter property
of $\Dis(X)$ can be easily derived from the following two lemmas,
first of which belongs to mathematical folklore.

\begin{lemma} A metric space $(M,d)$ is uniformly locally path
connected provided there is a dense subset $D\subset M$ such that
for any $\e>0$ there is $\delta>0$ such that any two points
$x,y\in D$ with $d(x,y)<\delta$ can be linked by a path
$f:[0,1]\to M$ with diameter $<\e$.
\end{lemma}

To apply this Lemma with $M=\Dis(X)$ and $D=\uDis(X)$ we need

\begin{lemma} Let $(X,d)$ be a uniformly locally path
equi-connected metric space having totally bounded small balls.
For any $\e>0$ there is $\delta>0$ such that any two sets
$A,B\in\uDis(X)$ with $d_H(A,B)<\delta$ can be linked by a path
$f:[0,1]\to\Dis(X)$ with diameter $<\e$.
\end{lemma}

\begin{proof} Since $X$ is uniformly locally path equi-connected,
there are $r >0$ and a continuity modulus $\w$ such that any two
points $x,y\in X$ with $d(x,y)<r $ can be connected by a path
$f:[0,d(x,y)]\to X$ such that $f(0)=x$, $f(d(x,y))=y$ and
$\w_f\le\w$.

We can assume that $r $ is so small that the ball $B(x,2\w(r ))$
is totally bounded for each $x\in X$.

Given arbitrary $\e>0$ find positive $\delta<r $ such that
$\w(\delta)<\e$. Fix any uniformly discrete sets $A,B\in\uDis(X)$
with $\rho=d_H(A,B)<\delta$.

For each points $a\in A$ select a point $b_a\in B$ with
$d(a,b_a)\le \rho$ and for any $b\in B$ find $a_b\in A$ with
$d(b,a_b)\le\rho$. By the choice of $r $ and $\w$, there are maps
$f_a,f_b:[0,\rho]\to X$ such that $f_a(0)=a$, $f_a(\rho)=b_a$,
$f_b(0)=a_b$, $f_b(\rho)=b$, and $\w_{f_a}\le\w$, $\w_{f_b}\le\w$.
Define a map $f:[0,\rho]\to \Dis(X)$ letting
$f(t)=\{f_a(t),f_b(t):a\in A,\; b\in B\}$ for $t\in[0,\rho]$. Let
us show that $f$ is well-defined, i.e, for each $t\in[0,\rho]$ the
set $f(t)$ is closed and discrete in $X$.

Given any point $x_0\in X$ it suffices to find a neighborhood of
$x_0$ having finite intersection with $f(t)$. For this consider
the sets $A'=\{a\in A:d(f_a(t),x_0)<\w(r )\}$ and $B'=\{b\in
B:d(f_b(t),x_0)<\w(r )\}$. Since $d(f_a(t),a)\le\w(\rho)\le\w(r )$
for all $a\in A$, the set $A'$ lies in the ball $B(x_0,2\w(r ))$.
Now the uniform discreteness of $A\supset A'$ and the total
boundedness of $B(x_0,2\w(r ))$ imply that the set $A'$ is finite.
By analogy, we can prove that the set $B'$ is finite.
Consequently, $B(x_0,\w(r ))\cap f(t)$ has size $\le|A'|+|B'|$,
which just yields that $f(t)$ is closed and discrete.

It is easy to see that the path $f:[0,\rho]\to \Dis(X)$ links the
sets $A$ and $B$ and satisfies $\w_f\le\w$. This yields also that
$\diam f([0,\rho])\le\w(\rho)\le\w(\delta)<\e$.
\end{proof}

\section{Local chain connectedness in hyperspaces}

In this section we will prove the implications (6)$\Rightarrow$(7)
of Theorem~\ref{t1} and $(4)\Rightarrow(5)$ of Theorem~\ref{t2}.
Our proof will exploit the famous Ramsey Theorem asserting that
for each finite coloring $\chi:[\IN]^2\to\{0,\dots,r\}$ of the
collection $[\IN]^2$ of two-element subsets of $\IN$ there is an
infinite subset $I\subset\IN$ such that the set $[I]^2$ is
monochromatic, that is $\chi([I]^2)\equiv {\mathrm const}$, see
\cite[\S1.5]{GRS}.

A sequence $(x_n)_{n=1}^\infty$ of points of a set $X$ will be
called {\em regular} if either $|\{x_n:n\in\IN\}|=1$ or else
$x_n\ne x_m$ for distinct $n,m\in\IN$.

\begin{lemma}\label{l1} A metric space $(X,d)$ fails to be
uniformly locally chain equi-connected at a subset $B\subset X$ if
and only if $\exists\e>0\;\forall \delta>0\;\exists\eta>0\;
\exists \{x_i,y_i:i\in\IN\}\subset B$ such that
\begin{enumerate}
\item for all $i\in\IN$ \;
$d(x_i,y_i)<\delta$ and $x_i,y_i$ cannot be
connected by an $\eta$-chain in $X$ of
length $\le i$ and diameter $<\e$;
\item the sets $\{x_i:i\in\IN\}$, $\{y_i:i\in\IN\}$ are
$\eta$-separated;
\item the sequences $(x_i)$ and $(y_i)$ are regular.
\end{enumerate}
\end{lemma}

\begin{proof} The ``if'' part is trivial. To prove the ``only if''
part, assume that $X$ fails to be uniformly locally chain
equi-connected at $B$. This means that for some $\e_0>0$ and each
$\delta>0$ there is $\eta>0$ such that for each $n\in\IN$ there is
a pair of points $a_n,b_n\in B$ with $d(a_n,b_n)<\delta$ that
cannot be linked by an $\eta$-chain of length $\le n$ and diameter
$\le\e_0$. Without loss of generality we can assume that
$\eta<\e_0/4$. Let $\e=\e_0/2$.

Define a coloring of $[\IN]^2$ into four colors as follows.
Colorate a two-element subset $\{i,j\}\subset\IN$ in
\begin{itemize}
\item {\em black} if $d(a_i,a_j)<\eta$ and
$d(b_i,b_j)<\eta$;
\item {\em white} if $d(a_i,a_j)\ge\eta$ and
$d(b_i,b_j)\ge\eta$;
\item {\em green} if $d(a_i,a_j)<\eta$ and
$d(b_i,b_j)\ge\eta$;
\item {\em blue} if $d(a_i,a_j)\ge\eta$ and
$d(b_i,b_j)<\eta$.
\end{itemize}
Applying the Ramsey Theorem, we can find an infinite subset
$N\subset [2,\infty)\cap\IN$ such that all two-element subsets of
$N$ are colored by the same color $\chi\in\{$black, white, green,
blue$\}$. Let $N=\{n_k:k\in\IN\}$ be an increasing enumeration of
$N$ and note that $n_k>k$ for all $k$. Fix any number $m\in N$.

1. Assuming that the color $\chi$ is black, let $x_i=a_{m}$,
$y_i=b_{m}$ for all $i\in N$. Then the conditions (2), (3) of
Lemma~\ref{l1} is trivially satisfied.

To verify the condition (1), assume that for some $i\in\IN$ the
points $x_i=a_{m}$ and $y_i=b_{m}$ can be connected by a
$\eta$-chain of diameter $<\e=\e_0/2$ and length $\le i$. Take any
$n\in N$ with $n>i+2$. Since
$\max\{d(a_{m},a_n),d(b_{m},b_n)\}<\eta<\e_0/4$, the points $a_n$
and $b_n$ can be connected by an $\eta$-chain of length $\le
i+2\le n$ and diameter $<\frac{\e_0}2+2\frac{\e_0}4=\e_0$, which
contradicts to the choice of these points.

2. Assuming that the color $\chi$ is white let $x_i=a_{n_i}$ and
$y_i=b_{n_i}$ for all $i\in \IN$ and observe that the sequences
$(x_i)$, $(y_i)$ satisfy the requirements of the lemma.

3. Assuming that the color $\chi$ is green let $x_i=a_{m}$ and
$y_i=b_{n_i}$ for all $i\in \IN$. It is clear that the sequences
$(x_i)$, $(y_i)$ satisfy the condition (2) and (3). Assume that
for some $i\in\IN$ the points $x_i=a_m$ and $y_i=b_{n_i}$ can be
connected by an $\eta$-chain of length $\le i$ and diameter
$\le\e_0/2$. Since $d(a_m,a_{n_i})<\eta<\e_0/2$, the points
$a_{n_i},b_{n_i}$ can be connected by an $\eta$-chain of length
$\le i+1\le n_i$ and diameter $\le\frac{\e_0}2+\eta\le\e_0$, which
contradicts to the choice of these points.

4. Assuming that the color $\chi$ is blue let $x_i=a_{n_i}$ and
$y_i=b_{m}$ for all $i\in \IN$. By analogy with the previous case
we can show that the sequences $(x_i)$ and $(y_i)$ satisfy the
requirements of the lemma.
\end{proof}

\begin{lemma}\label{lem7} A metric space $(X,d)$
is uniformly locally chain equi-connected at a closed subset
$B\subset X$, provided the hyperspace $2^X$ is locally chain
connected at $\Dis_0(B)$.
\end{lemma}

\begin{proof}
Assume that some closed subset $B$ of $X$ fails to be uniformly
locally chain equi-connected. Applying Lemma~\ref{l1}, for some
$\e_0>0$ and each $n\in\IN$ we can find $\eta_n>0$ and regular
sequences $(x_{n,k})_{k\in\IN}$, $(y_{n,k})_{k\in\IN}$ of points
in $B$ such that $d(x_{n,k},y_{n,k})<\frac1n$, the sets
$\{x_{n,k}:k\in\IN\}$ and $\{y_{n,k}:k\in\IN\}$ are
$\eta_n$-separated and for each $k\in\IN$ the points
$x_{n,k},y_{n,k}$ cannot be connected by an $\eta_n$-chain of
length $\le k$ and diameter $\le\e_0$. Without loss of generality,
we may assume that $\eta_{n}<\min\{\eta_{n-1},\frac1{n}\}$.

For every $n\in\IN$ find an infinite subset $I_n\subset \IN$ such
that $$\diam\{x_{n,k},y_{n,k}:k\in I_n\}<\eta_n+\inf\{
\diam\{x_{n,k},y_{n,k}:k\in I\}: I\subset \IN,\; |I|=\infty\}.$$
Let $A_n=\{x_{n,k},y_{n,k}:k\in I_n\}$, $n\in\IN$. Now consider
separately two cases.
\medskip

I. $\inf_{n\in\IN}\diam A_n>0$. In this case $$\Delta=\inf\{
\diam\{x_{n,k},y_{n,k}:k\in I\}: I\subset \IN,\; |I|=\infty\}>0$$
and for any subset $D\subset X$ of diameter $<\Delta$ the set
$\{k\in\IN:\{x_{n,k},y_{n,k}\}\subset D\}$ is finite for each $n$.
The condition $\lim_{n\to\infty}\eta_n=0<\inf_{n\in\IN}\diam A_n$
and the regularity of the sequences $(x_{n,k})_{k\in\IN}$,
$(y_{n,k})_{k\in\IN}$ imply that the sets $A_n$ are infinite for
all $n$ exceeding some number $n_0$.

It follows that for any finite subset $F\subset X$ and any $n>n_0$
the complement $A_n\setminus B(F,\Delta/3)$ is infinite. Using
this observation we shall construct a $\Delta/3$-separated subset
$D_0\subset\bigcup_{n\in\IN}A_n\subset B$ such that the
intersection $D_0\cap A_n$ is infinite for each $n>n_0$. The
construction of $D_0$ looks as follows. Let $\xi:\IN\to
(n_0,\infty)\cap\IN$ be a function such that $\xi^{-1}(n)$ is
infinite for all $n>n_0$. By induction for every $n\in\IN$ select
a point $d_n\in A_{\xi(n)}\setminus
B(\{d_1,\dots,d_{n-1},\Delta/3\})$. Then the set
$D_0=\{d_n:n\in\IN\}$ has the desired property.

Since $2^X$ is locally chain connected at $\Dis_0(B)\ni D_0$, for
the number $\e=\min\{\frac \Delta9,\frac{\e_0}{2}\}$ there is
$\delta>0$ such that each countable uniformly discrete subset
$D\subset B$ with $d_H(D_0,D)<\delta$ for any $\eta>0$ can be
connected with $D_0$ by an $\eta$-chain of diameter $<\e$. Take
any number $n>n_0$ with $\frac 1n<\min\{\frac \Delta9,\delta\}$
and consider two countable uniformly discrete sets: $$
\begin{aligned} D_x=&\;D_0\cup \{x_{n,k}\in A_n:
\{x_{n,k},y_{n,k}\}\cap D_0\ne\emptyset\}\\
 D_y=&\;D_0\cup
\{y_{n,k}\in A_n: \{x_{n,k},y_{n,k}\}\cap D_0\ne\emptyset\}
\end{aligned}
$$ It follows from $d(x_{n,k},y_{n,k})<\frac1n<\delta$ that
$d_H(D_0,D_x)\le\frac1n<\delta$ and consequently, $D_0$ can be
linked with $D_x$ by an $\eta_n$-chain $D_0,D_1,\dots,D_p=D_x$ of
diameter $<\e$. By analogy, $D_0$ can be connected with $D_y$ by
an $\eta_n$-chain $D_0=E_0,E_1,\dots,E_q=D_y$ of diameter $<\e$.

Find a number $l>p+q$ such that $\{x_{n,l},y_{n,l}\}\cap D_0$
contains some point $z$ (such a number $l$ exists because the
intersection $D_0\cap A_n$ is infinite). Since $x_{n,l}\in D_x$,
selecting suitable points in the $\eta_n$-chain
$D_0,D_1,\dots,D_p$ we can construct an $\eta_n$-chain
$x_{n,l}=z_0,\dots,z_p$ with $z_i\in D_{p-i}$ for $i\le p$. We
claim that $z_p=z$ and $\{z_0,\dots,z_p\}\subset B(z,\e_0/2)$.
Indeed, $\{z_0,\dots,z_p\}\subset D_0\cup\dots\cup D_p\subset
B(D_0,\e)\subset B(D_0,\e_0/2)$. Taking into account that $\e\le
\frac \Delta9$, the set $D_0$ is $\Delta/3$-separated,
$d(z_0,z)<\frac1n<\frac \Delta9$, and
$d(z_{i-1},z_i)<\eta_n<\frac1n<\frac \Delta9$ for $i\le p$, we
conclude that $\{z_0,\dots,z_p\}\subset B(z,\e)$ and thus $z_p\in
D_0\cap B(z,\e)=\{z\}$ which just yields $z_p=z$.

By analogy, using the $\eta_n$-chain $E_0,\dots,E_q$ between $D_0$
and $D_y$, we can construct an $\eta_n$-chain
$z_p,z_{p+1},\dots,z_{p+q}\in B(z,\e)$ connecting $z$ with
$y_{n,l}$. Then $z_0,\dots,z_{p+q}\in B(z,\e)$ is an
$\eta_n$-chain of length $p+q<l$ and diameter $<2\e\le \e_0$
linking the points $x_{n,l}$ and $y_{n,l}$, which contradicts to
their choice. This completes the proof in the case $\Delta>0$.
\bigskip

II. Now we assume that $\inf_{n\in\IN} \diam A_n=0$. Take an
infinite subset $N\subset\IN$ such that $\lim_{n\in N}\diam
A_n=0$. In each set $A_n$, $n\in N$, fix a point $a_n\in A_n$. We
consider separately three subcases.
\medskip

II.1. The set $\{a_n:n\in N\}$ is not totally bounded in $X$.
Then, replacing $N$ by a smaller subset, if necessary, we can
assume that the set $D_0=\{a_n:n\in N\}$ is $\alpha$-separated for
some $\alpha>0$. Since $2^X$ is locally chain connected at
$\Dis_0(B)\ni D_0$, for the number $\e=\min\{\frac \alpha
3,\frac{\e_0}{2}\}$ there is $\delta>0$ such that each countable
uniformly discrete subset $D\subset B$ with $d_H(D_0,D)<\delta$
can be connected with $D_0$ in $2^X$ by an $\eta$-chain of
diameter $<\e$ for every $\eta>0$. Take any $n\in N$ with $\diam
A_n<\min\{\frac \alpha 3,\delta\}$ and consider two countable
uniformly discrete sets: $$D_x=D_0\cup \{x_{n,i}: i\in I_n\}\mbox{
and } D_y=D_0\cup\{y_{n,i}: i\in I_n\}.$$ It follows from $\diam
A_n<\delta$ that $d_H(D_0,D_x)<\delta$ and consequently, $D_0$ can
be linked with $D_x$ and $D_y$ by $\eta_n$-chains of diameter
$<\e$ and length $\le l$ for some $l$. Using these chains and
repeating the argument from the case (I), for each $i\in I_n$ we
can construct an $\eta_n$-chain of diameter $<\e_0$ and length
$\le 2l$ connecting the points $x_{n,i},y_{n,i}$, which is not
possible for $i>2l$.
\medskip

II.2. Next we consider the case when the sequence $\{a_n:n\in N\}$
has a cluster point $a_\infty\in B$. Replacing $N$ by a smaller
subset we can assume that the sequence $\{a_n:n\in N\}$ converges
to $a_\infty$. Consider the one-point subset $D_0=\{a_\infty\}$
and use the local chain connectedness of $2^X$ at $\Dis_0(B)$ in
$2^X$ to find $\delta>0$ such that any countable uniformly
discrete subset $D\subset B$ with $d_H(D_0,D)<\delta$
 can be linked with $D_0$ by an $\eta$-chain of
diameter $\e_0/2$ for each $\eta>0$. Find $n\in N$ such that
$d(a_n,a_\infty)<\delta/2$ and $\diam A_n<\delta/2$. Consider two
countable bounded uniformly discrete sets: $$D_x=\{x_{n,i}: i\in
I_n\}\mbox{ and } D_y=\{y_{n,i}: i\in I_n\}.$$ It follows that
$d_H(D_0,D_x)\le d(a_n,a_\infty)+\diam A_n<\delta$ and
consequently, $D_0=\{a_\infty\}$ can be linked with $D_x$ and
$D_y$ by $\eta_n$-chains in $2^X$ of diameter $<\e_0/2$ and length
$\le l$ for some $l$. Using these chains for each $i\in I_n$ we
can construct an $\eta_n$-chain in $B(a_\infty,\e_0/2)$ of length
$\le 2l$ connecting the points $x_{n,i},y_{n,i}$, which is not
possible for $i>2l$.
\medskip

II.3. Finally consider the case of a totally bounded sequence
$\{a_n:n\in N\}$ having no limit points in $X$. Replacing $N$ by a
smaller subset, if necessary, we may assume that the sequence
$\{a_n:n\in N\}$ is Cauchy and has diameter $<\e_0/8$. With help
of the Ramsey theorem we select from this sequence an especially
nice subsequence as follows. Colorate a two-point subset
$\{n,m\}\in [N]^2$ with $n<m$ in
\begin{itemize}
\item {\em white} if $a_n$ and $a_m$ can be connected by an
$\eta_m$-chain of diameter $<\e_0/4$;
\item {\em black\/} otherwise.
\end{itemize}
Apply the Ramsey Theorem to find an infinite subset $M\subset N$
such that all two-element subsets of $M$ have the same color. We
claim that this color is white. For this it suffices to find a
two-element subset $\{n,m\}\subset M$ colored in white.

Consider the totally bounded closed discrete subset
$D_0=\{a_n:n\in M\}\in \Dis_0(B)$. Since $2^X$ is locally chain
connected at $\Dis_0(B)$, for the real number $\frac{\e_0}{16}$
there is $\delta>0$ such that each finite subset $D\subset B$ with
$d_H(D_0,D)<\delta$ can be connected with $D_0$ by an $\eta$-chain
of diameter $<\frac{\e_0}{16}$ for each $\eta>0$. Since
$D_0=\{a_n:n\in M\}$ is totally bounded, there is a finite subset
$F\subset M$ such that $d_H(D_0,D)<\delta$ where $D=\{a_n:n\in
F\}$. Take any $m\in M$ with $m>\max F$ and find an $\eta_m$-chain
of diameter $<\frac{\e_0}{16}$ connecting $D_0$ and $D$. Using
this chain construct an $\eta_m$-chain connecting $a_m\in D_0$
with some point $a_n\in D$, $n\in F$, and lying in the
$\frac{\e_0}{16}$-ball around $D_0$. This means that the diameter
of the latter chain is $<\diam{D_0}+2\frac{\e_0}{16}\le
\frac{\e_0}4$ and thus the set $\{n,m\}$ is white.

For the real number $\delta$ defined above find a number $n\in
M\setminus F$ such that $\diam A_n<\delta$. For this number $n$
find a finite subset $E\subset M$ such that
$d_H(D_0,D)<\min\{\eta_n,\delta/2\}$ where $D=\{a_i:i\in E\}$. By
the choice of $M$ for any numbers $i<j$ in $E$ the points
$a_i,a_j$ are connected by an $\eta_j$-chain of diameter
$<\e_0/4$. Let $L$ be the maximal length of these chains.

Consider two countable uniformly discrete subsets of $B$:
$$D_x=D\cup\{x_{n,k}: \in I_n\}\mbox{ and } D_y=D\cup\{y_{n,k}:
k\in I_n\}.$$ It follows that
$\max\{d_H(D_0,D),d_H(D_0,D_x),d_H(D_0,D_y)\}<\delta$ and
consequently, $D$ can be linked with $D_x$ and $D_y$ by
$\eta_n$-chains in $B(D_0,\frac{\e_0}{16})$ with length $\le l$
for some $l$. Take any number $k>2(L+l)$ and using these chains
construct an $\eta_n$-chain of length $\le l$ connecting the point
$x_{n,k}$ with some point $x\in D$ and an $\eta_n$-chain of length
$\le l$ connecting the point $y_{n,k}$ with some point $y\in D$.
These chains lie in $B(D_0,\frac{\e_0}{16})$. Taking into account
that the point $a_n\in D_0$ can be connected with the points
$x,y\in D$ by $\eta_n$-chains of length $\le L$ and diameter
$<\e_0/4$ we conclude that the points $x_{n,k}$ and $y_{n,k}$ can
be connected by an $\eta_n$-chain of length $2(L+l)<k$, which lies
in $B(D_0,{\e_0}4)$ and hence has diameter $<\diam
D_0+2\frac{\e_0}4\le\frac58\e_0<\e_0$, which contradicts to the
choice of the points $x_{n,k},y_{n,k}$.
\end{proof}

To prove the implication (8)$\Rightarrow$(4) of Theorem~\ref{t1}
we will exploit the following simple lemma allowing us to
transform traces into paths in hyperspaces. Below by $\cl_X(A)$ we
denote the closure of a subset $A\subset X$.

\begin{lemma}\label{lem8} Let $X$ be a metric space and $f:Q\to 2^X$ be a
uniformly continuous map defined on a countable dense subset
$Q\supset\{0,a\}$ of $[0,a]$. Then the map $g:[0,a]\to 2^X$
defined by $g:t\mapsto\cl_X(\cup\{f(q):|q-t|\le\dist(t,\{0,a\})$
is a continuous path such that $g(0)=f(0)$, $g(a)=f(a)$ and
$\w_g\le \w_f$.
\end{lemma}

\section{Proof of theorem~\ref{t2}}

Let $(X,d)$ be a metric space and $\mathcal H\subset 2^X$ be an
open subspace such that $\Dis_0(F)\subset\mathcal H$ for all
$F\in\mathcal H$. The implication $(2)\Rightarrow(1)$ follows from
Proposition~\ref{prop1} while
$(1)\Rightarrow(2)\Rightarrow(3)\Rightarrow(4)$ are trivial. The
implications $(4)\Rightarrow(5)$ and $(5)\Rightarrow(6)$ follow
from Lemmas~\ref{lem7} and \ref{lem2}, respectively.

Now we verify the implication $(6)\Rightarrow(2)$. To show that
$\mathcal H$ is locally path-connected, fix any element
$F\in\mathcal H$ and $\e>0$. Since the set $\mathcal H$ is open in
$2^X$ there is $r>0$ such that any $A\in 2^X$ with $d_H(A,F)<2r$
belongs to $\mathcal H$. In particular, the closed
$r$-neighborhood $\overline{B}(F,r)=\{x\in X:d(x,F)\le r\}$
belongs to $\mathcal H$ and thus $X$ is uniformly locally trace
equi-connected at $\overline{B}(F,r)$ by (6). Consequently, there
are $\e_0>0$ and a continuity modulus $\w$ such that for any
points $x,y\in \overline{B}(F,r)$ with $d(x,y)<\e_0$ and any
countable dense subset $Q\supset\{0,d(x,y)\}$ of $[0,d(x,y)]$
there is a uniformly continuous function $f:Q\to X$ with $f(0)=x$,
$f(d(x,y))=y$ and $\w_f\le \w$. Without loss of generality we may
assume that $\w(\e_0)<\min\{\e,r\}$.

We claim that any element $A\in\mathcal H$ with $d_H(A,F)<\e_0$
can be linked with $F$ by a path $g:[0,\e_0]\to \mathcal H$ lying
in the $\e$-ball around $F$ in $2^X$. Pick any countable dense
subset $Q\supset\{0,\e_0\}$ of $[0,\e_0]$. For any points $a\in
A$, $b\in F$ find points $b_a\in F$, $a_b\in A$ with
$d(a,b_a)<\e_0$ and $d(b,a_b)<\e_0$. For the obtained pairs
$(a,b_a)$, $(b,a_b)$ fix uniformly continuous maps $f_a:Q\to X$,
$f_b:Q\to X$ such that $f_a(0)=a$, $f_a(\e_0)=b_a$, $f_b(0)=a_b$,
$f_b(\e_0)=b$, and $\w_{f_a}\le\w$, $\w_{f_b}\le\w$. Define a
trace $f:Q\to 2^X$ letting $f(q)=\cl_X(\{f_a(q),f_b(q):a\in
A,\;b\in B\})$ for $q\in Q$. It can be shown that $f(0)=A$,
$f(\e_0)=F$, and $\w_f\le\w$. Using Lemma~\ref{lem8} produce a
continuous path $g:[0,\e_0]\to 2^X$ defined by
$g(t)=\cl_X(\cup\{f(q):|q-t|\le\dist(t,\{0,\e_0\})\}$ for
$t\in[0,\e_0]$. For this path we get $g(0)=A$, $g(\e_0)=F$, and
$\w_g\le\w_f\le\w$, which implies that $g([0,\e_0])$ lies in the
closed $\w(\e_0)$-ball around $F$ in $2^X$. Since
$\w(\e_0)<\min\{\e,r\}$ this ball lies in $\mathcal H$, which
means that $\mathcal H$ is locally path connected.

\section{Proof of Theorem~\ref{t1}}

The implication $(4)\Rightarrow(3)$ follows from
Proposition~\ref{prop2} while
$(3)\Rightarrow(1)\Rightarrow(2)\Rightarrow(5)\Rightarrow(6)$ are
trivial. The implication $(6)\Rightarrow(7)$ is proven in
Lemma~\ref{lem7} while
$(7)\Rightarrow(8)\Rightarrow(9)\Rightarrow(7)$ in
Lemmas~\ref{lem2}, \ref{lem5} and \ref{lem3}, respectively. The
implication $(8)\Rightarrow(4)$ can be proved by analogy with the
proof of $(6)\Rightarrow(2)$ from Theorem~\ref{t2}. For a complete
$X$ the equivalence $(8)\Leftrightarrow(10)$ is trivial.

The implication $(11)\Rightarrow(3)$ is trivial. Assuming that $X$
is uniformly locally compact we shall prove the implication
$(10)\Rightarrow(11)$. Assuming that $X$ is uniformly locally path
equi-connected, and applying Proposition~\ref{prop3}, we conclude
that the space $\Dis(X)$ is a uniform ANR. According to a result
of K.Sakai \cite{Sa}, a metric space is uANR provided it contains
a dense uANR. Consequently, each space $\mathcal H\subset 2^X$
containing $\Dis(X)$ is uANR.

Finally, assuming that $X$ is chain equi-connected, we shall
verify the implications
$(3)\Rightarrow(13)\Rightarrow(12)\Rightarrow(1)$. Assume that
$2^X$ is a uANR. Since the space $X$ is chain equi-connected we
can apply Proposition 4.6 of \cite{CK} to conclude that $2^X$ is
chain connected (in should be mentioned that in \cite{CK} chain
equi-connected sets are called uniformly $C$-connected). Now
Proposition~\ref{prop2} implies that $2^X$ is a uniform AR. This
completes the proof of $(3)\Rightarrow(13)$. Two other
implications are trivial.

\section{Proof of Theorem~\ref{t4}}

Assume that for each bounded subset $B$ of $X$ the space $X$ is
uniformly locally chain equi-connected at $B$ and $B$ lies in a
bounded chain equi-connected subspace of $X$. By
Corollary~\ref{t3}, the hyperspace $\Bd(X)$ is an ANR and $X$ is
uniformly locally trace equi-connected at each bounded subset of
$X$. In light of Proposition~\ref{prop1} to show that $\Bd(X)$ is
an AR, it suffices to verify that the space $\Bd(X)$ is path
connected. For this we shall show that each one-point subset
$\{x_0\}\in 2^X$ can linked with a non-empty closed bounded subset
$A\in 2^X$ by a continuous path in $\Bd(X)$. Our assumption
implies that $\{x_0\}\cup A$ lies in a bounded chain
equi-connected subspace $C\subset X$.

Since $X$ is uniformly locally trace equi-connected at the bounded
set $C$ there is $\e_0<1$ and a continuity modulus $\w$ such that
any two points $x,y\in C$ with $d(x,y)<\e_0$ can be linked by a
trace $f:\IQ\cap[0,1]\to X$ defined on the space of rational
numbers of $[0,1]$ such that $f(0)=x$, $f(1)=y$, and $\w_f\le \w$.

Use the chain equi-connectedness of $C$ to find $l\in\IN$ such
that any points $x,y\in C$ can be linked by an $\e_0$-chain of
length $\le l$. Hence for each point $a\in A\subset C$ we can pick
a map $f_a:\{0,\dots,l\}\to C$ such that $f_a(0)=x_0$, $f_a(l)=a$,
and $d(f_a(i-1),f_a(i))<\e_0$. By the choice of $\e_0$ and $\w$ we
can extend $f_a$ to a uniformly continuous map
$f_a:\IQ\cap[0,l]\to X$ whose continuity modulus does not exceed
$l\w$. Then the function $f:\IQ\cap[0,l]\to \Bd(X)$ defined by
$f:q\mapsto \cl_X(\cup\{f_a(q):a\in A\})$ is a trace in $\Bd(X)$
linking $\{x_0\}$ with $A$. Finally, use Lemma~\ref{lem8} to
transform $f$ into a continuous path in $\Bd(X)$ connecting
$\{x_0\}$ and $A$. Therefore $\Bd(X)$ is a path connected ANR.
Applying Proposition~\ref{prop1} we conclude that $\Bd(X)$ is an
AR.

Now assume that $\Bd(X)$ is an AR. By Corollary~\ref{t3}, $X$ is
uniformly locally chain equi-connected at each bounded subset $B$.
Since $\Bd(X)$ is path connected we can apply Proposition 4.7 of
\cite{CK} to conclude that each bounded subset of $X$ lies in a
bounded chain equi-connected subspace of $X$ (in should be
mentioned that in \cite{CK} chain equi-connected sets are called
uniformly $C$-connected).

\section{Acknowledgements}

The authors would like to express their sincere thanks to Katsuro
Sakai for valuable remarks concerning the initial version of this
manuscript.

\end{document}